# On the gradient of Schwarz symmetrization of functions in Sobolev spaces

Marco Bramanti

*Dipartimento di Matematica, Politecnico di Milano. Via Bonardi 9. 20133 Milano. Italy.*

**Sunto.** Sia $\mathcal{S}$ uno spazio di Sobolev o Orlicz-Sobolev di funzioni non necessariamente nulle al bordo del dominio. Si danno condizioni sufficienti su una funzione non negativa in $\mathcal{S}$ affinché la sua simmetrizzata di Schwarz appartenga ancora ad $\mathcal{S}$. Questi risultati sono ottenuti per mezzo di disuguaglianze isoperimetriche relative e generalizzano in un certo senso un noto teorema di Polya-Szegö. Si dimostra anche che il riarrangiamento di una qualsiasi funzione in $\mathcal{S}$ è localmente in $\mathcal{S}$.

**Abstract.** Let $\mathcal{S}$ be a Sobolev or Orlicz-Sobolev space of functions not necessarily vanishing at the boundary of the domain. We give sufficient conditions on a nonnegative function in $\mathcal{S}$ in order that its spherical rearrangement ("Schwartz symmetrization") still belongs to $\mathcal{S}$. These results are obtained via relative isoperimetric inequalities and somewhat generalize a well-known Polya-Szegö's theorem. We also prove that the rearrangement of any function in $\mathcal{S}$ is locally in $\mathcal{S}$.

If u is a nonnegative function in $H^{1,2}(\Re^n)$, u has compact support, and $\tilde{u}$ denotes the Schwarz symmetrization of u, then a well known theorem by Polya-Szegö states that $\tilde{u}$ belongs to $H^{1,2}(\Re^n)$ and:

$$\int |D\tilde{u}|^2 dx \leq \int |Du|^2 dx. \qquad (*)$$

(Henceforth, we will indicate with D the gradient of a function of n variables or the derivative of a function of one real variable).

In particular, this formula holds for $u \in H_0^{1,2}(\Omega)$, where $\Omega$ is a bounded domain of $\Re^n$, the first integral is taken on the ball $\tilde{\Omega}$ having the same measure of $\Omega$ and the second is taken on $\Omega$.

If u is a function in $H^{1,2}(\Omega)$, not necessarily vanishing at the boundary, or if u belongs to $H_0^{1,2}(\Omega)$ but assumes also negative values (and so does $\tilde{u}$), then inequality (*) can actually fail, and $\tilde{u}$ does not necessarily belong to $H^{1,2}(\Omega)$ (see examples below). So, a natural question is under which additional assumptions a nonnegative function in $H^{1,2}(\Omega) \setminus H_0^{1,2}(\Omega)$ has Schwarz symmetrization in $H^{1,2}(\tilde{\Omega})$. In section 1 we will prove some different sufficient conditions (in terms of the size of the set on which u vanishes) in order to a Polya-Szegö-type estimate holds, that is:

$$\int_{\tilde{\Omega}} |D\tilde{u}|^2 dx \leq (\text{const.}) \cdot \int_{\Omega} |Du|^2 dx.$$

Moreover, we will prove that whenever u is an $H^{1,2}(\Omega)$ function (even of changing sign), $\tilde{u}$ belongs to $H_{\text{loc}}^{1,2}(\tilde{\Omega})$ and for any ball $\tilde{\Omega}_\epsilon$ concentric to $\tilde{\Omega}$ and with measure $|\Omega| - \epsilon$, one has:

$$\int_{\tilde{\Omega}_\epsilon} |D\tilde{u}|^2 dx \leq c(\epsilon) \cdot \int_{\Omega} |Du|^2 dx.$$

where c does not depend on u. (See section 2). All these results can naturally be generalized to Orlicz-Sobolev spaces. This will be done in section 3.

The interest in studying properties of the rearrangement of functions in $H^{1,2}(\Omega)$, or vanishing on part of the boundary, comes from the application of symmetrization techniques to elliptic or parabolic P.D.E. with boundary conditions of Neumann or mixed type: so thm. 2.1 and corollary 2.2 have been used in investigating parabolic Neumann problems, see [2]. We also mention [8], in which a similar result to thm. 1.3 is stated, in a different context: this result is related to the study of elliptic mixed problems, which is carried out in [13].

### Some notations and examples

If u is a real measurable function defined on $\Omega$, we define:
the distribution function of u:

$$\mu(t) = \left|\left\{x \in \Omega : u(x) > t\right\}\right| \quad \text{for } t \in \Re \qquad (0.1)$$

( $|\ |$ denotes Lebesgue measure);

the decreasing rearrangement of u:

$$u^*(s) = \inf\{t \in \Re: \mu(t) \leq s\} \text{ for } s \in [0, |\Omega|]; \qquad (0.2)$$

the Schwarz symmetrization of u:

$$\tilde{u}(x) = u^*(c_n |x|^n) \text{ for } x \in \tilde{\Omega}, \qquad (0.3)$$

where $\tilde{\Omega}$ is the sphere centred at the origin with the same measure of $\Omega$; $c_n$ is the measure of the unit ball in $\Re^n$.

For general properties of these functions, see [12]; note that, in our definition, $u^*$ and $\tilde{u}$ assume also negative values, if u is a function of changing sign, whereas rearrangements are sometimes defined for $|u|$.

From (0.3) it follows:

$$|D\tilde{u}(x)| = nc_n |Du^*(c_n|x|^n)| \cdot |x|^{n-1}$$
$$\int_{\tilde{\Omega}} |D\tilde{u}(x)|^2 dx = (nc_n^{1/n})^2 \int_0^{|\Omega|} |Du^*(s)|^2 s^{2-2/n} ds. \qquad (0.4)$$

Hence, if $\tilde{u} \in H^{1,2}(\tilde{\Omega})$, $u^* \in H^{1,2}(\epsilon, |\Omega|)$ for any $\epsilon > 0$, so that $u^* \in AC(\epsilon, |\Omega|)$ for any $\epsilon > 0$.

For better understanding the problem of assuring integrability of $|D\tilde{u}|^2$, let us consider the case of a radially symmetric and *increasing* function u defined on a ball $\Omega$, i.e.:

$$u(x) = u^*(|\Omega| - c_n |x|^n). \qquad (0.5)$$

In this case one has:

$$\int_\Omega |Du(x)|^2 dx = (nc_n^{1/n})^2 \int_0^{|\Omega|} |Du^*(s)|^2 (|\Omega| - s)^{2-2/n} ds. \qquad (0.6)$$

Comparing (0.4) and (0.6) one sees how it may happen that $u \in H^{1,2}(\Omega)$ but $\tilde{u} \notin H^{1,2}(\tilde{\Omega})$. Take, for instance, $u^*(s) = \sqrt{|\Omega| - s}$ and u as in (0.5). Then:

$$\int_\Omega |Du(x)|^2 dx = \frac{(nc_n^{1/n})^2}{4} \int_0^{|\Omega|} s^{1-2/n} ds < \infty \text{ for every } n \geq 2, \text{ while:}$$

$$\int_{\tilde{\Omega}} |D\tilde{u}(x)|^2 dx = \frac{(nc_n^{1/n})^2}{4} \int_0^{|\Omega|} \frac{s^{2-2/n}}{|\Omega|-s} ds = \infty \text{ for every } n.$$

Similarly, if one defines: $u^*(s) = \sqrt{|\Omega| - s} - \sqrt{|\Omega|}$ and u as in (0.5), one has an example of a (negative) function $u \in H_0^{1,2}(\Omega)$ such that $\tilde{u} \notin H^{1,2}(\tilde{\Omega})$.

**Remark 0.1.** The above example works for $n \geq 2$. If $n = 1$ inequality (*) can actually be proved for any nonnegative function in $H^{1,2}(\Omega)$. (See [6], p.35). So in this paper we will always consider $n \geq 2$.

## 1. Isoperimetric inequalities and $\mathcal{L}^2$ norm of the gradient of $\tilde{u}$

Here we want to obtain a proof of integrability of $|D\tilde{u}|^2$ without assuming that u vanishes at the boundary of $\Omega$. In what follows u will be a *nonnegative* function defined on $\Omega$. A first basic tool we need is Federer's "coarea formula", as appears in [11]:

if $f \in \mathcal{L}^1(\Re^n)$ and v is a nonnegative Lipschitz function with compact support, then:

$$\int_{\Re^n} f(x) |Dv(x)| dx = \int_0^{+\infty} dt \int_{\{x: v(x)=t\}} f(x) dH_{n-1}(x). \qquad (1.1)$$

(Here and below, $H_{n-1}$ stands for $(n-1)$-dimensional Hausdorff measure).

Let us consider a nonnegative Lipschitz function u defined on $\Omega$. If $\Omega$ is Lipschitz, we can extend u to a compact supported Lipschitz function on $\Re^n$. Then, if $f \in \mathcal{L}^1(\Omega)$ and we put $f \equiv 0$ outside $\Omega$, (1.1) becomes:

$$\int_\Omega f(x) |Du(x)| dx = \int_0^{+\infty} dt \int_{\{x \in \Omega: u(x)=t\}} f(x) dH_{n-1}(x). \qquad (1.2)$$

From (1.2) it follows in particular:

$$\int_{\{x \in \Omega: u(x)>t\}} |Du(x)| dx = \int_t^{+\infty} H_{n-1}\{x \in \Omega: u(x) = \xi\} d\xi. \qquad (1.3)$$

Note that:

$$\left\{x \in \Omega: u(x) = \xi\right\} \supseteq \partial\left\{x \in \Omega: u(x) > \xi\right\} \cap \Omega, \text{ and:} \tag{1.4}$$

$$H_{n-1}\left\{x \in \Omega: u(x) = \xi\right\} \geq P_\Omega\left\{x \in \Omega: u(x) > \xi\right\}.$$

Here $P_\Omega$ stands for the perimeter, in the sense of De Giorgi, relative to $\Omega$. For a definition of this concept in the general case, see [9]. However, we will only use the fact that $P_\Omega(E) \leq H_{n-1}(\partial E) \cap \Omega$ for every measurable subset E of $0\Omega$, and, if $\partial E$ is sufficiently smooth, this is an equality. (See [4]). The perimeter of E, P(E), is equal to $P_\Omega(E)$ when $\Omega = \Re^n$. We recall De Giorgi's isoperimetric inequality in $\Re^n$:

$$P(E) \geq n\, c_n^{1/n}\, |E|^{1-1/n}.$$

The next theorem points out the role of isoperimetric inequalities in Polya-Szegö-type estimates.

**Theorem 1.1.** Let $\Omega$ be a bounded Lipschitz domain in $\Re^n$, $n \geq 2$, $u \in \text{Lip}(\Omega)$, $u \geq 0$ in $\Omega$, and assume that u satisfies:

$$P_\Omega\left\{x \in \Omega: u(x) > t\right\} \geq \gamma \cdot \mu(t)^{1-1/n} \tag{1.5}$$

for some positive constant $\gamma$, any $t \geq 0$. (Here and below, $\mu$ is the distribution function of u, defined in (0.1)).
Then $\tilde{u} \in \text{Lip}(\tilde{\Omega})$, and:

$$\int_{\tilde{\Omega}} |D\tilde{u}|^2\, dx \leq \left(\frac{n\, c_n^{1/n}}{\gamma}\right)^2 \int_\Omega |Du|^2\, dx. \tag{1.6}$$

**Proof.** (Here we revise an argument of [11]). Let us prove that $\tilde{u}$ is Lipschitz. If L is a constant such that $|Du(x)| \leq L$ in $\Omega$, and t, h such that $0 < h < t$, then:

$$L[\mu(t-h) - \mu(t)] \geq \int_{\{x \in \Omega:\, t-h < u(x) \leq t\}} |Du(x)|\, dx \quad = \quad \text{(by (1.3), (1.4))}$$

$$= \int_{t-h}^{t} P_\Omega\left\{x \in \Omega: u(x) > \xi\right\} d\xi \geq \text{ (by (1.5)) } \gamma \int_{t-h}^{t} \mu(\xi)^{1-1/n}\, d\xi \geq$$

$$\geq \text{ (by monotonicity of } \mu\text{) } \gamma \cdot h \cdot \mu(t)^{1-1/n}.$$

Hence $\mu$ is strictly decreasing in $(0, \|u\|_\infty)$, so that $u^*$ is continuous and satisfies:

$$u^*(s) - u^*(s+k) \leq \frac{L}{\gamma}\, s^{-1+1/n} \cdot k$$

for any $k > 0$, $s, s+k \in (0, |\Omega|)$. Therefore $u^* \in \text{AC}(\epsilon, |\Omega|)$ for any $\epsilon > 0$ and:

$$0 \leq -\frac{du^*}{ds}(s) \leq \frac{L}{\gamma}\, s^{-1+1/n}. \tag{1.7}$$

By the definition of $\tilde{u}$ and (1.7) one can compute:

$$|\tilde{u}(x) - \tilde{u}(y)| = \left|\int_{c_n|y|^n}^{c_n|x|^n} \frac{du^*}{ds}(s)\, ds\right| \leq L\, \frac{n\, c_n^{1/n}}{\gamma}\, |y - x|$$

that is $\tilde{u}$ is Lipschitz in $\Omega$.

Let us prove now that (1.6) holds. From (1.3)-(1.4) it follows:
$$-\frac{d}{dt} \int_{\{x \in \Omega:\, u(x) > t\}} |Du(x)|\, dx = P_\Omega\left\{x \in \Omega: u(x) > t\right\} \geq \text{ (by (1.5)) } \gamma \cdot \mu(t)^{1-1/n}. \tag{1.8}$$

From (1.2) it follows that:

$$\varphi(t) \equiv \int_{\{x \in \Omega:\, u(x) > t\}} |Du(x)|^2\, dx = \int_{t}^{+\infty} d\xi \int_{\{x \in \Omega:\, u(x) = \xi\}} |Du|\, dH_{n-1}(x)$$

from which one reads that $\varphi$ is absolutely continuous, so that:

$$\int_\Omega |Du|^2\, dx = \varphi(0) = \int_0^{+\infty} -\varphi'(t)\, dt. \tag{1.9}$$

Writing differential quotients and applying Holder's inequality one has:

$$-\varphi\prime(t) \geq \frac{1}{-\mu\prime(t)}\left[-\frac{d}{dt}\int_{\{x\in\Omega:\,u(x)>t\}} |Du(x)|\,dx\right]^2. \qquad (1.10)$$

From (1.8), (1.9), (1.10) it follows:

$$\int_\Omega |Du|^2\,dx \geq \gamma^2 \int_0^{+\infty} \frac{\mu(t)^{2-2/n}}{-\mu\prime(t)}\,dt. \qquad (1.11)$$

Now consider $\tilde{u}$. Since its level sets are balls, in (1.10) the equal sign holds, and (1.8) becomes:

$$-\frac{d}{dt}\int_{\{x\in\tilde{\Omega}:\,\tilde{u}(x)>t\}} |D\tilde{u}|\,dx \;=\; n\,c_n^{1/n}\,\mu(t)^{1-1/n}$$

Hence:

$$\int_{\tilde{\Omega}} |D\tilde{u}|^2\,dx \;=\; \left(n\,c_n^{1/n}\right)^2 \int_0^{+\infty} \frac{\mu(t)^{2-2/n}}{-\mu\prime(t)}\,dt. \qquad (1.12)$$

From (1.11)-(1.12) it follows estimate (1.6). □

Now we are interested in discussing sufficient conditions in order that (1.5) holds. In the following the function u is still supposed nonnegative and Lipschitz in $\overline{\Omega}$.

(*i*)   If $u = 0$ on $\partial\Omega$, we obtain Polya-Szegö's theorem, since:

$$P_\Omega\{x\in\Omega:\,u(x)>t\} \;=\; P\{x:\,u(x)>t\} \;\geq\; n\,c_n^{1/n}\cdot\mu(t)^{1-1/n},$$

by the isoperimetric inequality in $\Re^n$. So $\gamma = n\,c_n^{1/n}$, and (1.6) holds with constant equal to 1.

(*ii*)   Suppose that: $|\text{support of } u| \leq \frac{|\Omega|}{2}$. The *relative* isoperimetric inequality of $\Omega$ says that:

$$Q\cdot P_\Omega(E) \;\geq\; \min\Big(|E|,\,|\Omega\setminus E|\Big)^{1-1/n} \qquad (1.13)$$

for some constant $Q > 0$, any measurable set $E \subseteq \Omega$. (Such an inequality certainly holds if $\Omega$ is Lipschitz). Then:

$$P_\Omega\{x\in\Omega:\,u(x)>t\} \;\geq\; Q^{-1}\mu(t)^{1-1/n}, \qquad (1.14)$$

and (1.5) holds with $\gamma = Q^{-1}$.

(*iii*)   More generally, suppose that:

$$\Big|\{x\in\Omega:\,u(x)=0\}\Big| \;=\; \epsilon$$

with $0 < \epsilon < \frac{|\Omega|}{2}$. Fix $t > 0$. If $\mu(t) \leq \frac{|\Omega|}{2}$, (1.14) still holds. Otherwise, from (1.13) we get:

$$Q\cdot P_\Omega\{x\in\Omega:\,u(x)>t\} \;\geq\; \Big(|\Omega|-\mu(t)\Big)^{1-1/n} \;\geq\; [\alpha\,\mu(t)]^{1-1/n}$$

with $\alpha = \frac{\epsilon}{|\Omega|-\epsilon}$. Hence (1.5) holds with:

$$\gamma \;=\; Q^{-1}\alpha^{1-1/n} \;=\; Q^{-1}\cdot\left(\frac{\epsilon}{|\Omega|-\epsilon}\right)^{1-1/n}$$

and (1.6) holds with constant:

$$\left(\frac{Q\,n\,c_n^{1/n}}{\alpha^{1-1/n}}\right)^2.$$

(*iv*)   Now, suppose that:

$$H_{n-1}\Big\{x\in\partial\Omega:\,u(x)=0\Big\} \;=\; \epsilon > 0.$$

We also suppose that $\Omega$ satisfies the following geometric property (this already appears in [10]):

$$H_{n-1}\Big(\partial E \cap \partial\Omega\Big) \;\leq\; \mathcal{C}\cdot P_\Omega(E) \qquad (1.15)$$

for some positive constant $\mathcal{C}$, for any measurable $E \subseteq \Omega$ such that $|E| \leq \frac{|\Omega|}{2}$. (If $\Omega$ is Lipschitz, (1.15) actually holds).

Fix $t > 0$. Again, we consider the case $\mu(t) > \frac{|\Omega|}{2}$; then, by (1.15):

$$H_{n-1}\left(\partial\{x \in \Omega: u(x) \leq t\} \cap \partial\Omega\right) \leq \mathcal{C} \cdot P_\Omega\{x \in \Omega: u(x) \leq t\}.$$

Hence:

$$P_\Omega\{x \in \Omega: u(x) > t\} = P_\Omega\{x \in \Omega: u(x) \leq t\} \geq$$

$$\geq \tfrac{1}{\mathcal{C}} H_{n-1}\left(\partial\{x \in \Omega: u(x) \leq t\} \cap \partial\Omega\right) \geq \tfrac{1}{\mathcal{C}} H_{n-1}\{x \in \partial\Omega: u(x) = 0\} =$$

$$= \tfrac{\epsilon}{\mathcal{C}} \geq \tfrac{\epsilon}{\mathcal{C}} \left(\tfrac{\mu(t)}{|\Omega|}\right)^{1-1/n}.$$

So (1.5) holds with:

$$\gamma = \min\left(Q^{-1}, \tfrac{\epsilon}{\mathcal{C}|\Omega|^{1-1/n}}\right)$$

and (1.6) holds with constant:

$$\max\left\{\left(Q\,n\,c_n^{1/n}\right), \left(\tfrac{\mathcal{C}\,n\,c_n^{1/n}|\Omega|^{1-1/n}}{\epsilon}\right)\right\}^2.$$

Note that: $\tfrac{n c_n^{1/n}|\Omega|^{1-1/n}}{\epsilon} \leq \tfrac{H_{n-1}(\partial\Omega)}{\epsilon},$ which is a more expressive ratio.

(v) Suppose that $E = \{x \in \Omega: u(x) = 0\}$ is such that its projection on at least one hyperplane has positive $(n-1)$-dimensional Hausdorff measure, in symbols:

$$H_{n-1}\left(\Pi(E)\right) = \epsilon > 0 \text{ for some projection } \Pi.$$

For any $t > 0$, the set $A = \{u \leq t\}$ contains $E$, so:

$$P(A) = H_{n-1}(\partial A) \geq H_{n-1}\left(\Pi(A)\right) \geq H_{n-1}\left(\Pi(E)\right) = \epsilon.$$

Now, if $\mu(t) > \tfrac{|\Omega|}{2}$, one has:

$$\mathcal{C} \cdot P_\Omega\{u > t\} \geq H_{n-1}\left(\partial\{u \leq t\} \cap \partial\Omega\right).$$

Hence:

$$P_\Omega\{u > t\} \geq \tfrac{1}{\mathcal{C}+1} P\{u \leq t\} \geq \tfrac{\epsilon}{\mathcal{C}+1} \geq \tfrac{\epsilon}{\mathcal{C}+1} \tfrac{\mu(t)^{1-1/n}}{|\Omega|^{1-1/n}}.$$

So (1.5) holds with:

$$\gamma = \min\left(Q^{-1}, \tfrac{\epsilon}{(\mathcal{C}+1)|\Omega|^{1-1/n}}\right).$$

Now we state separately the results obtained from (iii)-(iv)-(v).

**Theorem 1.2.** Let $\Omega$ be a bounded Lipschitz domain in $\Re^n$, $n \geq 2$; let $u \in H^{1,2}(\Omega)$, $u \geq 0$ in $\Omega$, and suppose that $|\text{support of } u| = |\Omega| - \epsilon$ for some $\epsilon > 0$. Then $\tilde{u} \in H_0^{1,2}(\tilde{\Omega})$ and:

$$\int_{\tilde{\Omega}} |D\tilde{u}|^2 dx \leq L^2 \cdot \int_\Omega |Du|^2 dx \tag{1.16}$$

with $L = \left(\tfrac{Q\,n\,c_n^{1/n}}{\alpha^{1-1/n}}\right)$, where $Q$ is as in (1.13) and: $\alpha = \tfrac{\epsilon}{|\Omega|-\epsilon}$ if $\epsilon \leq \tfrac{|\Omega|}{2}$; $\alpha = 1$ otherwise.

**Theorem 1.3.** Let $\Omega$ be as above, let $u \geq 0$, $u \in \mathcal{H}$, where $\mathcal{H}$ is the closure in $H^{1,2}$-norm of the space:

$$\mathcal{H}\prime \equiv \left\{\varphi \in \text{Lip}(\Omega): \text{supp}\,\varphi \cap F = \emptyset\right\}$$

and F is a fixed closed subset of $\partial\Omega$ with $H_{n-1}(F) = \epsilon > 0$. Then $\tilde{u} \in H^{1,2}_0(\tilde{\Omega})$ and (1.16) holds with:

$$L = \max\left\{\left(Q\,n\,c_n^{1/n}\right), \left(\tfrac{\mathcal{C}\,n\,c_n^{1/n}|\Omega|^{1-1/n}}{\epsilon}\right)\right\} \tag{1.17}$$

and $\mathcal{C}$ as in (1.15).

**Theorem 1.4.** Let $\Omega$ be as above, let $u \geq 0$, $u \in \mathcal{H}$, where $\mathcal{H}$ is the closure in $H^{1,2}$-norm of the space:

$$\mathcal{H}\prime \equiv \left\{\varphi \in \text{Lip}(\Omega): \text{supp}\,\varphi \subseteq \overline{\Omega}\backslash F\right\}$$

where F is a closed subset of $\Omega$ with the property stated in (v). Then $\tilde{u} \in H^{1,2}_0(\tilde{\Omega})$ and (1.16) holds with:

$$L = \max\left\{\left(Q\,n\,c_n^{1/n}\right), \left(\tfrac{(\mathcal{C}+1)\,n\,c_n^{1/n}|\Omega|^{1-1/n}}{\epsilon}\right)\right\}.$$

**Remark 1.5.** We note that the spaces $\mathcal{H}$ defined in thms. 1.3-1.4 are properly contained in $H^{1,2}(\Omega)$ whenever F has positive capacity. This is the case, in particular, if F has positive $(n-1)$-measure. Moreover, if F has (positive and) finite $(n-1)$-measure and is a *regular* set in the sense of geometric measure theory (that is a. e. $(H_{n-1})$ point of F is a density point in sense $H_{n-1}$) then property (v) is certainly satisfied. (See [5], p.87).

**Proof of theorem 1.2.** If $u \in H^{1,2}(\Omega)$, $u \geq 0$ and $\Omega$ is Lipschitz, u may be approximated in $H^{1,2}$-norm with smooth functions $u_m$ in $\overline{\Omega}$. (See [1], thm. 3.18). Moreover, if the support of u has measure $|\Omega| - \epsilon$, then for any $\epsilon_1 \in (0,\epsilon)$ $\{u_m\}$ can be choosen such that:

$$\left|\left\{x \in \Omega: u_m(x) = 0\right\}\right| \geq \epsilon_1.$$

Hence, for every m, $u_m$ satisfies (1.16) (with $\epsilon$ replaced by $\epsilon_1$), so that $\{\tilde{u}_m\}$ is a bounded sequence in $H^{1,2}_0(\tilde{\Omega})$. Let $\tilde{u}_m$ be a subsequence converging to some $v \in H^{1,2}_0(\tilde{\Omega})$ weakly in $H^{1,2}$ and strongly in $\mathcal{L}^2$. By [3], $u_m \to u$ in $\mathcal{L}^2(\Omega)$ implies $\tilde{u}_m \to \tilde{u}$ in $\mathcal{L}^2(\tilde{\Omega})$, so $v \equiv \tilde{u}$ and $\tilde{u} \in H^{1,2}_0(\tilde{\Omega})$. Then from weak convergence it follows that u satisfies (1.16) for any $\epsilon_1 < \epsilon$, and hence for $\epsilon$, too. $\square$

**Proof of theorem 1.3.** If $u_m \in \mathcal{H}\prime$, $u_m \to u$ in $H^{1,2}(\Omega)$, then $u_m$ satisfies (1.16)-(1.17). Hence arguing as above, it follows that these hold for u. Note that the condition $\text{supp}\,u_m \cap F = \emptyset$ implies that $|\text{supp}\,u_m| < |\Omega|$; hence $\tilde{u}_m \in H^{1,2}_0(\tilde{\Omega})$, and so does u. $\square$

In a similar way it follows theorem 1.4. Incidentally, we note that a Sobolev embedding theorem for functions vanishing on part of the boundary can be derived from thm. 1.3:

**Corollary 1.6.** Let $u \in \mathcal{H}$, where $\mathcal{H}$ is as in theorem 1.3 or 1.4. Then the following estimate holds:

$$\|u\|_{\mathcal{L}^{2^*}(\Omega)} \leq \text{const.} \cdot \|Du\|_{\mathcal{L}^2(\Omega)}. \tag{1.18}$$

**Proof.** It is sufficient to prove (1.18) for $u \geq 0$. Then $\tilde{u} \in H^{1,2}_0(\tilde{\Omega})$, so by Sobolev's embedding theorem and theorem 1.3 (or 1.4) one has:

$$\|u\|_{\mathcal{L}^{2^*}(\Omega)} = \|\tilde{u}\|_{\mathcal{L}^{2^*}(\tilde{\Omega})} \leq C \cdot \|D\tilde{u}\|_{\mathcal{L}^2(\tilde{\Omega})} \leq C \cdot \|Du\|_{\mathcal{L}^2(\Omega)} \quad \square$$

## 2. Local integrability of $|D\tilde{u}|^2$ for $u \in H^{1,2}(\Omega)$

Theorem 1.2 allows us to prove the following result, which holds for *any* function $u \in H^{1,2}(\Omega)$ (even assuming negative values):

**Theorem 2.1.** Let $\Omega$ be as in theorem 1.2, $u \in H^{1,2}(\Omega)$. Then $\tilde{u} \in H^{1,2}_{\text{loc}}(\tilde{\Omega})$ and, for any $\epsilon > 0$, one has:

$$\int_{\tilde{\Omega}_\epsilon} |D\tilde{u}|^2 dx \leq c(\epsilon) \cdot \left(Q n c_n^{1/n}\right)^2 \cdot \int_\Omega |Du|^2 dx$$

where $\tilde{\Omega}_\epsilon$ is the sphere centred at the origin with measure $|\Omega| - \epsilon$ and:

$$c(\epsilon) = \left(\frac{|\Omega|-\epsilon}{\epsilon}\right)^{2-2/n} \quad \text{if } \epsilon \leq \frac{|\Omega|}{2}, \ c(\epsilon) = 1 \text{ otherwise.}$$

Moreover, $u^* \in AC(\epsilon, |\Omega| - \epsilon)$.

**Proof.** Put $h = u^*(\frac{|\Omega|}{2})$, and let $u_1$, $u_2$ be the positive and negative parts of $(u - h)$. Then $u_i \in H^{1,2}(\Omega)$, $|\text{supp } u_i| \leq \frac{|\Omega|}{2}$ ($i = 1,2$). So by theorem 1.2 $\tilde{u}_i \in H_0^{1,2}(\tilde{\Omega})$ and:

$$\int_{\tilde{\Omega}} |D\tilde{u}_i|^2 dx \leq \left(Q n c_n^{1/n}\right)^2 \cdot \int_\Omega |Du_i|^2 dx.$$

In particular, $u_i^* \in AC(\epsilon, |\Omega|)$ for any $\epsilon > 0$. Now, noting that:

$$(v^+)^*(s) = (v^*)^+(s) \tag{2.1}$$

$$(v^-)^*(s) = (v^*)^-(|\Omega| - s) \tag{2.2}$$

one has:

$$(u^* - h)^+ \in AC(\epsilon, |\Omega|), \ (u^* - h)^- \in AC(0, |\Omega| - \epsilon), \text{ so that:}$$

$$u^* \in AC(\epsilon, |\Omega| - \epsilon) \text{ for any } \epsilon > 0. \tag{2.3}$$

Note also that:

$$\widetilde{(u - h)}^+ = \widetilde{(u - h)^+} \tag{2.4}$$

whereas the same is *not* true for the *negative* part. To handle the gradient of $\widetilde{(u - h)}^-$, let us observe that, for any $\epsilon > 0$, one has, by (0.4):

$$\int_{\tilde{\Omega}_\epsilon} |D\widetilde{(u-h)}^-|^2 dx = \left(n c_n^{1/n}\right)^2 \int_{\frac{1}{2}|\Omega|}^{|\Omega|-\epsilon} s^{2-2/n} |Du^*(s)|^2 ds \tag{2.5}$$

while, by (0.4) and (2.2):

$$\int_{\tilde{\Omega}_\epsilon} |D\widetilde{(u-h)^-}|^2 dx = \left(n c_n^{1/n}\right)^2 \int_\epsilon^{\frac{1}{2}|\Omega|} s^{2-2/n} |Du^*(|\Omega| - s)|^2 ds = \tag{2.6}$$

$$= \left(n c_n^{1/n}\right)^2 \int_{\frac{1}{2}|\Omega|}^{|\Omega|-\epsilon} (|\Omega| - s)^{2-2/n} |Du^*(s)|^2 ds.$$

Comparing (2.5) and (2.6) we can write:

$$\int_{\tilde{\Omega}_\epsilon} |D\widetilde{(u-h)}^-|^2 dx \leq \left(\frac{|\Omega|-\epsilon}{\epsilon}\right)^{2-2/n} \int_{\tilde{\Omega}_\epsilon} |D\widetilde{(u-h)^-}|^2 dx. \tag{2.7}$$

Finally, we can estimate:

$$\int_{\tilde{\Omega}_\epsilon} |D\tilde{u}(x)|^2 dx = \int_{\tilde{\Omega}_\epsilon} |D\widetilde{(u-h)}^+|^2 dx + \int_{\tilde{\Omega}_\epsilon} |D\widetilde{(u-h)}^-|^2 dx \leq \text{ (by (2.4), (2.7))}$$

$$\leq \int_{\tilde{\Omega}} |D\tilde{u}_1|^2 dx + \left(\frac{|\Omega|-\epsilon}{\epsilon}\right)^{2-2/n} \int_{\tilde{\Omega}} |D\tilde{u}_2|^2 dx \leq \text{ (by (2.1))}$$

$$\leq \left(Q n c_n^{1/n}\right)^2 \cdot \max\left(1, \left(\tfrac{|\Omega|-\epsilon}{\epsilon}\right)^{2-2/n}\right) \cdot \left\{\int_\Omega |Du_1|^2 dx + \int_\Omega |Du_2|^2 dx\right\} =$$

$$= \left(Q n c_n^{1/n}\right)^2 \cdot \max\left(1, \left(\tfrac{|\Omega|-\epsilon}{\epsilon}\right)^{2-2/n}\right) \cdot \int_\Omega |Du|^2 dx.$$

So the theorem is completely proved. □

From the previous theorem it follows the next estimate, giving an approximation result for rearrangements:

**Corollary 2.2.** Let $\Omega$ be as above, u, v $\in H^{1,2}(\Omega)$. Then for any $\epsilon > 0$ one has:

$$\sup_{s \in (\epsilon, |\Omega| - \epsilon)} |(u^* - v^*)(s)| \leq c_1(n, Q, |\Omega|) \|u - v\|_2 +$$

$$+ c_2(\epsilon, n, Q, |\Omega|) \|u - v\|_2^{1/2} \cdot \left\{ \|Du\|_2 + \|Dv\|_2 \right\}^{1/2}.$$

In particular, if $u_m$ is a sequence of $H^{1,2}$ functions converging to u in $H^{1,2}(\Omega)$, then $u_m^*$ converges to $u^*$ uniformly in $(\epsilon, |\Omega| - \epsilon)$ for any $\epsilon > 0$.

**Proof.** We start by noting that if $\varphi$ is an absolutely continuous function on [a,b], then:

$$\varphi(s) \leq \frac{1}{b-a} \int_a^b \varphi(\sigma) \, d\sigma + \int_a^b |\varphi'(\sigma)| \, d\sigma$$

for every $s \in [a,b]$. Applying this formula to the function:

$$\varphi(s) = \left(u_m^*(s) - u^*(s)\right)^2 \quad \text{on the interval } (\epsilon, |\Omega| - \epsilon)$$

we have, by Hölder's inequality:

$$\left(u_m^*(s) - u^*(s)\right)^2 \leq \frac{1}{|\Omega| - 2\epsilon} \int_\epsilon^{|\Omega|-\epsilon} |u_m^*(\sigma) - u^*(\sigma)|^2 \, d\sigma +$$

$$+ 2 \left( \int_\epsilon^{|\Omega|-\epsilon} [\sigma^{-1+1/n} \, |u_m^*(\sigma) - u^*(\sigma)|]^2 \, d\sigma \right)^{\frac{1}{2}} \cdot \left( \int_\epsilon^{|\Omega|-\epsilon} [\sigma^{1-1/n} \, |\frac{du_m^*}{d\sigma}(\sigma) - \frac{du^*}{d\sigma}(\sigma)|]^2 \, d\sigma \right)^{\frac{1}{2}} \equiv$$

$$\equiv A_m + 2 B_m \cdot C_m. \tag{2.8}$$

Now:

$$A_m \leq \frac{1}{|\Omega| - 2\epsilon} \|u_m - u\|_2^2 \quad \text{and:} \tag{2.9}$$

$$B_m \leq \epsilon^{-1+1/n} \|u_m - u\|_2, \quad \text{while:} \tag{2.10}$$

$$C_m \leq \left( \int_0^{|\Omega|-\epsilon} [\sigma^{1-1/n} \, |\frac{du_m^*}{d\sigma}(\sigma)|]^2 \, d\sigma \right)^{\frac{1}{2}} + \left( \int_0^{|\Omega|-\epsilon} [\sigma^{1-1/n} \, |\frac{du^*}{d\sigma}(\sigma)|]^2 \, d\sigma \right)^{\frac{1}{2}} =$$

$$= \left( \int_{\tilde{\Omega}_\epsilon} |D\tilde{u}_m|^2 \, dx \right)^{\frac{1}{2}} + \left( \int_{\tilde{\Omega}_\epsilon} |D\tilde{u}|^2 \, dx \right)^{\frac{1}{2}} \leq \text{ (by thm. 2.1)}$$

$$\leq c(n, Q, |\Omega|) \cdot \left\{ \|Du_m\|_2 + \|Du\|_2 \right\}. \tag{2.11}$$

Collecting (2.8), (2.9), (2.10), (2.11) one gets the result. □

## 3. Extension to Orlicz-Sobolev spaces

Let A: $[0, +\infty) \to [0, +\infty)$ be an "N-function" (see [7]), that is A is an increasing continuous convex function, such that:

$$\lim_{t \to 0} \frac{A(t)}{t} = 0; \quad \lim_{t \to +\infty} \frac{A(t)}{t} = +\infty.$$

By Jensen's inequality, we can repeat the proof of theorem 1.1 and obtain, under the same assumptions:

$$\int_\Omega A(|Du(x)|) \, dx \geq \int_0^{+\infty} A\left( \frac{\gamma \mu(t)^{1-1/n}}{-\mu'(t)} \right) (-\mu'(t)) \, dt \tag{3.1}$$

$$\int_{\tilde{\Omega}} A(|D\tilde{u}(x)|) \, dx = \int_0^{+\infty} A\left( \frac{n c_n^{1/n} \mu(t)^{1-1/n}}{-\mu'(t)} \right) (-\mu'(t)) \, dt \tag{3.2}$$

We can rewrite (3.1)-(3.2) replacing A(r) with $A(\frac{r}{\lambda})$ for any fixed $\lambda > 0$.

Then, choosing $\lambda_0 = \frac{n c_n^{1/n}}{\gamma}$ we get:

$$\int_{\widetilde{\Omega}} A\left(\frac{|D\widetilde{u}|}{\lambda_0}\right) dx \leq \int_{\Omega} A(|Du|) dx. \tag{3.3}$$

Now, recall that the natural norm in the Orlicz space:

$$\mathcal{L}_A(\Omega) \equiv \left\{ u:\Omega \to \Omega, u \text{ measurable such that } \int_{\Omega} A\left(\frac{|u|}{\lambda}\right) dx < +\infty \text{ for some } \lambda > 0 \right\}$$

$$\text{is: } \|u\|_A \equiv \inf\left\{ \lambda > 0: \int_{\Omega} A\left(\frac{|u|}{\lambda}\right) dx \leq 1 \right\}.$$

Rewriting again (3.3) with $A(r)$ replaced by $A(\frac{r}{\lambda})$, and choosing $\lambda = \|Du\|_A$ we get:

$$\int_{\widetilde{\Omega}} A\left(\frac{|D\widetilde{u}(x)|}{\lambda_0 \|Du\|_A}\right) dx \leq 1. \text{ Hence:}$$

$$\|D\widetilde{u}\|_{\mathcal{L}_A(\widetilde{\Omega})} \leq \left(\frac{n c_n^{1/n}}{\gamma}\right) \|Du\|_{\mathcal{L}_A(\Omega)}. \tag{3.4}$$

So we have proved the following:

**Theorem 3.1.** Let $\Omega$ and $A$ be as above, let u be a nonnegative Lipschitz function in $\Omega$, such that one of the following holds:

(i)    $u = 0$ in $E \subseteq \Omega$ with $|E| = \epsilon > 0$
(ii)    $u = 0$ in $F \subseteq \partial\Omega$ with $H_{n-1}(F) = \epsilon > 0$
(iii)    $u = 0$ in $G \subseteq \Omega$ with $H_{n-1}\left(\Pi(G)\right) = \epsilon > 0$ for some projection $\Pi$ (see section 1).

Then $\widetilde{u} \in \text{Lip}(\widetilde{\Omega})$ and (3.4) holds, with $\gamma$ possibly depending on $n, \mathcal{C}, Q, |\Omega|, \epsilon$.

**Remark 3.2.** We did not state the previous theorem for $u \in H^1 \mathcal{L}_A(\Omega)$ because to apply a limit process as in the proof of theorems (1.2)-(1.3)-(1.4) we have to know that a bounded sequence in $H^1 \mathcal{L}_A(\Omega)$ has a weakly converging subsequence. This cannot be assured without further assumptions on A. To discuss this fact, we recall some results from the theory of Orlicz-Sobolev spaces. (See [1]).

We say that A satisfies a "global $\Delta_2$-condition" if:

$$A(2t) \leq \delta A(t) \text{ for some } \delta > 0, \text{ any } t > 0. \tag{3.5}$$

We say that A satisfies a "$\Delta_2$-condition near infinity" if (3.5) holds only for any $t \geq t_0$, for some $t_0 > 0$. We say that $(A,\Omega)$ is $\Delta$-regular if: A satisfies a global $\Delta_2$-condition, or: A satisfies a $\Delta_2$-condition near infinity and $|\Omega| < +\infty$. If $(A,\Omega)$ is $\Delta$-regular, then $\mathcal{L}_A(\Omega)$ and $H^1\mathcal{L}_A(\Omega)$ are reflexive spaces; if $\Omega$ is Lipschitz then $\mathcal{C}^{\infty}(\overline{\Omega})$ is dense in $H^1\mathcal{L}_A(\Omega)$; if $|\Omega| < +\infty$ then $\mathcal{L}_A(\Omega)$ is continuously embedded in $\mathcal{L}^1(\Omega)$. Using these facts one can repeat the proofs of theorems (1.2)-(1.3)-(1.4) to get the following:

**Theorem 3.3.** Let $\Omega$, A be as above. Suppose that A satisfies a $\Delta_2$-condition near infinity, and let u satisfy the assumptions of one of theorems 1.2, 1.3, 1.4, with $H^{1,2}(\Omega)$ replaced by $H^1\mathcal{L}_A(\Omega)$. Then $\widetilde{u} \in H^1_0 \mathcal{L}_A(\Omega)$ and (3.4) holds, with $\gamma$ possibly depending on $n, \mathcal{C}, Q, |\Omega|, \epsilon$.

**Example.** An example of Orlicz-Sobolev space which does not reduce to a standard Sobolev space and satisfies the previous theorem is the one defined by $A(r) = r^p \log(1+r)$ with $p \geq 1$.

Now we are interested in stating an analogue of theorem 2.1 for Orlicz-Sobolev spaces. We first consider the case of a Lipschitz function u. The analogue of formula (0.4) is:

$$\int_{\widetilde{\Omega}} A(|D\widetilde{u}(x)|) dx = \int_0^{|\Omega|} A\left(n c_n^{1/n} |Du^*(s)| s^{1-1/n}\right) ds.$$

Arguing as in section 2 one gets:

$$\int_{\widetilde{\Omega}_\epsilon} A\left(|D(u-h)^-|\right) dx = \int_{\frac{1}{2}|\Omega|}^{|\Omega|-\epsilon} A\left(n c_n^{1/n} |Du^*(s)| s^{1-1/n}\right) ds \tag{3.6}$$

$$\int_{\tilde{\Omega}_\epsilon} A\Big(\,|\,D\,(u-h)^-\sim\,|\,\Big)\,dx \;=\; \int_{\frac{1}{2}|\Omega|}^{|\Omega|-\epsilon} A\Big(n\,c_n^{1/n}\,|\,Du^*(s)\,|\,(\,|\,\Omega\,|\,-s)^{1-1/n}\Big)\,ds. \tag{3.7}$$

Comparing (3.6)-(3.7) one can write:

$$\int_{\tilde{\Omega}_\epsilon} A\Big(\tfrac{|D(u-h)\widetilde{\;\;}^-|}{\lambda}\Big)\,dx \;\leq\; \int_{\tilde{\Omega}_\epsilon} A\Big(\,|\,D\,(u-h)^-\sim\,|\,\Big)\,dx \tag{3.8}$$

with $\lambda = \left(\frac{|\Omega|-\epsilon}{\epsilon}\right)^{1-1/n}$ (we take $\epsilon < \frac{|\Omega|}{2}$, so $\lambda > 1$).

Applying (3.3) to the positive and negative parts of $(u-h)$ we get, by (3.8):

$$\int_{\tilde{\Omega}} A(\tfrac{|D\tilde{u}|}{\lambda})\,dx \;=\; \int_{\tilde{\Omega}_\epsilon} \Big\{A\Big(\tfrac{|D(u-h)\widetilde{\;\;}^+|}{\lambda}\Big) + A\Big(\tfrac{|D(u-h)\widetilde{\;\;}^-|}{\lambda}\Big)\Big\}\,dx \;\leq$$

$$\leq\; \int_{\tilde{\Omega}_\epsilon} \Big\{A\Big(\,|\,D\,(u-h)^+\sim\,|\,\Big) + A\Big(\,|\,D\,(u-h)^-\sim\,|\,\Big)\Big\}\,dx \;\leq$$

$$\leq\; \int_\Omega \Big\{A\Big(\lambda_0\,|\,D(u-h)^+\,|\,\Big) + A\Big(\lambda_0\,|\,D(u-h)^-\,|\,\Big)\Big\}\,dx \;=$$

$$=\; \int_\Omega A\Big(\lambda_0\,|\,Du\,|\,\Big)\,dx \quad \text{with } \lambda_0 = Q\,n\,c_n^{1/n}.$$

Again, rewriting the previous inequality for $A(\tfrac{r}{\rho})$ instead of $A(r)$ and choosing $\rho = \lambda_0 \|Du\|_A$ we find:

$$\|D\tilde{u}\|_{\mathcal{L}_A(\tilde{\Omega}_\epsilon)} \;\leq\; \Big(Q\,n\,c_n^{1/n}\Big)\left(\tfrac{|\Omega|-\epsilon}{\epsilon}\right)^{1-1/n} \|Du\|_{\mathcal{L}_A(\Omega)} \tag{3.9}$$

for every $\epsilon \in (0, \tfrac{|\Omega|}{2})$.

This holds for every Lipschitz function u defined in $\Omega$. From this fact we get, by approximation with smooth functions:

**Theorem 3.4.** Let $\Omega$, A be as in theorem 3.3. If $u \in H^1\mathcal{L}_A(\Omega)$ then $\tilde{u} \in H^1_{loc}\mathcal{L}_A(\Omega)$ and (3.9) holds. Moreover, $u^* \in AC_{loc}(\epsilon, |\Omega|-\epsilon)$ for every $\epsilon > 0$.